\documentclass[12pt]{amsart}
\usepackage{graphics,amssymb,amsfonts,amsmath}

\usepackage{mathrsfs}
\let\mathcal\mathscr
 
\def\lra{\longrightarrow}
 
\def\Z{{\bf Z}}

\def\P{{\bf P}}
\def\C{{\bf C}}

\def\N{{\bf N}}

\def\vide{\varnothing}

\def\Sing{\mathop{\rm Sing}\nolimits}
\def\base{\mathop{\rm Bs}\nolimits}
\def\div{\mathop{\rm div}\nolimits}
\def\Pic{\mathop{\rm Pic}\nolimits}

\def\codim{\mathop{\rm codim}\nolimits}

\def\Supp{\mathop{\rm Supp}\nolimits}

\def\ie{\hbox{i.e.}}
\def\isom{\simeq}
\def\iff{if and only if}

\def\ot{\otimes}

\def\av{abelian variety}
\def\avs{abelian varieties}
\def\ppav{principally polarized\ abelian\ variety}
\def\ppavs{principally polarized abelian\ varieties}
\def\pav{polarized abelian variety}
\def\pavs{polarized abelian varieties}

 \def\moins{\mathop{\hbox{\vrule height 3pt depth -2pt
width 5pt}\,}}

\def\PA{\Pic^0(A)}
\def\PB{\Pic^0(B)}
\def\PV{\Pic^0(V)}

\def\cF{{\mathcal F}}

\def\cI{{\mathcal I}}
\def\cJ{{\mathcal J}}
\def\cO{{\mathcal O}}
\def\a{{\alpha}}
\def\b{{\beta}}

\def\eps{{\varepsilon}}
\def\phi{{\varphi}}

\newtheorem{theo}{Theorem}
\newtheorem{prop}[theo]{Proposition}

\newtheorem{lemm}[theo]{Lemma}
\newtheorem{coro}[theo]{Corollary}

\newtheorem{rema}[theo]{Remark}

\begin{document} 

\title[Singularities of divisors on abelian varieties]{Singularities of divisors of low degree on abelian varieties}
\author{Olivier Debarre and Christopher Hacon}
\begin{abstract}
Building on previous work  of Koll\'ar, Ein, Lazarsfeld, and Hacon, we show that 
ample divisors of low degree on an abelian variety have mild singularities in 
case the abelian variety is simple or the degree of the polarization is two. 
 \end{abstract}
\address{Math\'ematique -- IRMA -- UMR 7501\\
ÊÊÊÊÊÊÊÊÊÊÊÊÊÊÊÊÊÊ Universit\'e Louis Pasteur\\
ÊÊÊÊÊÊÊÊÊÊÊÊÊÊÊÊÊÊ 7, rue Ren\'e Descartes\\
ÊÊÊÊÊÊÊÊÊÊÊÊÊÊÊÊÊÊ 67084 Strasbourg C\'edex, France}
\email{debarre@math.u-strasbg.fr}
\urladdr{http://www-irma.u-strasbg.fr/$\sim$debarre}

\address{Department of Mathematics \\ Ê
University of Utah\\ Ê
155 South 1400 East, Room 233\\
Salt Lake City, UT 84112, USA}
\email{hacon@math.utah.edu}
\urladdr{http://www.math.utah.edu/$\sim$hacon}

\thanks{Olivier Debarre was visiting the University of Michigan when part of this work was done, with support from William Fulton and Robert Lazarsfeld. Christopher Hacon was partially supported by NSA research grant 
no: ÊMDA904-03-1-0101 and by a grant from the Sloan Foundation.}
\maketitle

\vskip 1cm

\section{Introduction}

Since Koll\'ar used in \cite{kol}  the Kawamata--Viehweg
vanishing theorem to settle classical conjectures about
singularities of theta divisors in complex
\ppavs, the subject has known spectacular developments. Ein and Lazarsfeld, using
generic vanishing theorems of Green and Lazarsfeld, proved in
\cite{EL} that irreducible theta divisors are normal and gave an optimal bound
on   the dimension of the locus of points of given multiplicity of a multitheta
divisor. Hacon determined in \cite{hac1} exactly when this bound is attained
and obtained in \cite{hac2}  results for ample divisors of degree $2$.
 
In this article, we investigate more generally \avs\ with an indecomposable polarization of degree smaller than the dimension. When the degree increases, many special cases begin to appear, often due to the presence of reducible divisors that represent the polarization. In order to avoid overly technical statements, we restrict ourselves to two  cases: the case where the ambient \av\ is {\em simple,} and the case of polarizations of degree $2$ (thereby completing Hacon's above-mentioned results). Although we obtained almost complete results for polarizations of degree $3$, we   chose not to inflict their very technical proofs on the unsuspecting reader.

 We refer to Theorems
\ref{main},  
\ref{main2}, \ref{th6},  and  
\ref{th17}  for more precise formulations  and quote only the
following results. Let $(A,\ell)$ be a \pav\ of degree $d$ and dimension  $g$, and let $D$ be any effective divisor that represents $m\ell$, with $m>0$.

Assume $A$ is {\em simple}  and    $g>(d+1)^2/4 $.   {\em If $m=1$, the divisor $D$  is normal and has rational singularities; if $k\ge 2$,  the set of points of multiplicity at least $mk$ on
$D$  has codimension $>k$ in $A$.}

When $d=2$ and $(A,\ell)$ is only indecomposable, similar conclusions hold. Moreover, {\em  the set of points of multiplicity at least $mk$ on
$D$ has codimension $\le k$ in $A$
\iff\ 
$(A,\ell)$ is a double
\'etale cover of a product of at least
$k$  nonzero  \ppavs.}

The proofs of these results   systematically use generic vanishing theorems and precise descriptions of cohomological loci attached to various situations (see \S \ref{s2}). We work over the complex numbers.

\section{Singularities of pairs}

We just need a quick review of the basic terminology
relative to the singularities of a pair
$(A,D)$ consisting of an effective $\hbox{\bf Q}$-divisor $D$ in a smooth projective
variety $A$.   

A {\em log resolution} of the pair $(A,D)$ is a proper birational
morphism $\mu :A'\to A$ such that the union of
$\mu^{-1} (D )$ and the   exceptional locus of $\mu$
is a divisor with simple normal crossing support.
Write 
$$\mu^* (K_A+D )=K_{A'} +\sum a _i D _i$$
 where the Ê$D _i$ are distinct prime divisors on $A'$. 
The pair $(A,D )$ is
\begin{itemize}
\item {\em log canonical} if $a_i\le 1$ for all $i$;
\item {\em log terminal} if $a_i<1$ for all $i$;
\item {\em canonical} if $a_i\le 0$ for all $i$ such that $D_i$ is $\mu$-exceptional;
\end{itemize}
for some log resolution $\mu$.
The {\it multiplier ideal sheaf} associated to the pair $(A,D )$
is 
$$\cI (A,D)=\mu _*\bigl(\omega_{A'/A}(-\lfloor\mu ^*D \rfloor )\bigr)=
\mu _*\bigl(\cO _{A'}(-\sum \lfloor a _i\rfloor D _i)\bigr)$$
One sees that 
\begin{eqnarray*}
(A,D) \hbox{ log canonical}&\Longleftrightarrow& \cI(A,tD)=\cO_A \hbox{
 for all } t\in\hbox{\bf Q}\cap(0,1)
\\
(A,D)\hbox{ log terminal}&\Longleftrightarrow& \cI(A,D)=\cO_A
\end{eqnarray*}
Assume now that $D$ is a prime divisor
 in $A$.
The {\em adjoint ideal sheaf}
$\cJ(A,D)\subset\cO_A$ is defined in \cite{EL}, Proposition 3.1. For any
desingularization
$f:X\to D$, it fits into an exact sequence
\begin{equation}\label{adj}
0\to \omega_A\to \omega_A(D)\otimes\cJ(A,D) \to
f_*\omega_X\to 0
\end{equation}
 of sheaves on $A$ and (\cite{EL}, Proposition 3.1; \cite{kolsc},  Corollary
7.9.2, Theorem 7.9, and Theorem (11.1.1))
\begin{eqnarray*}
(A,D) \hbox{ canonical}&\Longleftrightarrow& \cJ(A,D)=\cO_A
\\
&\Longleftrightarrow& D\hbox{ is normal and has rational
singularities}\nonumber
\end{eqnarray*}
Furthermore, for any positive integers $m\ge1$ and $k\ge2$, we have
\begin{eqnarray}\label{sing}
(A, \textstyle\frac1m D) \hbox{ log canonical}&\Longrightarrow&
 \lfloor \textstyle\frac{1}{m+1} D\rfloor =0\ {\rm and}\nonumber\\
 && \codim_A(\Sing_{mk}D)\ge k\nonumber
\\
(A,\textstyle\frac1m D)\hbox{  log terminal} &\Longrightarrow&
 \lfloor \textstyle\frac1m D\rfloor =0 \ {\rm and}\\
 &&   \codim_A(\Sing_{mk}D)> k\nonumber\\
(A,D)\hbox{  canonical }&\Longrightarrow&
\codim_A(\Sing_kD)> k\nonumber
\end{eqnarray}

\section{Polarized abelian
varieties}

Let $A$ be an \av\ of dimension $g$. Any line bundle $L$ on  
$A$ induces a morphism
$\phi_L:A\to\PA$ defined  by
$\phi_L(a)=\tau_a^*L\ot L^{-1}$, where $\tau_a:A\to A$ is the translation
$x\mapsto x-a$. This morphism only depends on the numerical equivalence class $[L]$ of $L$ and
will also be denoted by $\phi_{[L]}$.
 We denote its kernel by
$K(L)$ or $K([L])$. If $D$ is a divisor on $A$, we write $[D]$ for $[\cO_A(D)]$ and $K(D)$ for $K([D])$. The line bundle $L$ is
 ample \iff\ $K(L)$ is finite, in which case this group has order $d^2$, where
$$d=h^0(A,L)=\frac{1}{g!}c_1(L)^g$$
 is the {\em degree} of $L$. A polarization on $A$ is a
numerical  equivalence class of ample line bundles on $A$.
A polarization of degree $1$ is
called {\em principal} and a divisor representing it is called a theta
divisor. A polarization $\ell$ is of type $(d)$ if
$K(\ell)\isom (\hbox{\bf Z}/d\hbox{\bf Z})^2$. If $d$ is {\em prime}, any polarization of degree
$d$ is of type $(d)$. 

A \pav\
$(A,\ell)$ is {\em indecomposable} if it is not the product of nonzero polarized \avs. If   $g\ge 2$,   a general element of
$\ell$ is prime.

\section{Singularities of ample divisors in  abelian
varieties}

This section contains the central results of this article. Some auxiliary results will be proved later in \S\ \ref{s2} and \S\ \ref{app}.

Let $(A,\ell)$ be a \pav. We  study the singularities of a divisor  in $m\ell$, with $m>0$,  when the dimension is large enough with respect to the degree.

\begin{theo}\label{main}
Let $(A,\ell)$ be a {\em simple} polarized abelian variety of degree $d$ and dimension $g> (d+1)^2/4$.
\begin{itemize}
\item[{\rm a)}] Every divisor in $\ell$ is prime, normal, and has
rational singularities.
\item[{\rm b)}]
If $m\ge2$ and $D$ is a divisor in $m\ell$, the pair
$(A,\frac1m D)$ is  log terminal unless $D=mE$, with $E\in \ell$.
\end{itemize}
\end{theo}

In particular, in case b), the pair $(A,\frac1m D)$ is  log canonical. In view of further investigations, it is natural to   conjecture that the conclusion of the theorem hold under the weaker assumption  $g>d$.
 We can prove the conjecture   for $d\le 3$. For $g>4=d $, we can prove that the pair
$(A,\frac1m D)$ is  log terminal when $A$ is general. For $5\ge g\ge d$, we can prove that the pair
$(A,\frac1m D)$ is  log canonical when $A$ is simple. However, because of the technical nature of our arguments, 
we do not pursue this here.

Note that, in any dimension $\ge 2$, and for any $d\ge 3$, there are   examples (obtained by the construction of Remark \ref{rem10}) of indecomposable (but not simple!) \pavs\ $(A,\ell)$ of   degree $d $ and, for any $m\ge d-1$, of pairs  $(A,\frac1m D)$ that are {\em not} log canonical because $D$ has a component of multiplicity $\bigl[ \frac{md}{d-1}\bigr]>m$ (see (\ref{sing})).

However, for polarizations of degree  $2$, the  results of Theorem \ref{main} can be extended to the case where the \pav\ is only indecomposable.

\begin{theo}[Degrees $1$ and $2$]\label{main2}
 Let $(A,\ell)$ be an {\em indecomposable}
 \pav\  of degree $d\le 2$ and dimension   $g >d$. 
\begin{itemize}
\item[{\rm a)}] Every prime divisor in $\ell$ is normal and has
rational singularities.
\item[{\rm b)}] If $m\ge2$ and $D$ is a divisor in $m\ell$ such that $\lfloor \frac1m D\rfloor =0$, the pair
$(A,\frac1m D)$ is  log terminal.
\end{itemize}
\end{theo}

Note that $\ell$ may very well contain reducible elements (see \S \ref{app}).  As to b), we explain in   Corollary \ref{ccc} exactly when   the assumption $\lfloor \frac1m D\rfloor =0$ fails to hold (recall that in any event,   the pair $(A,\frac1m D)$
is always log canonical when $g\ge d$, as proved in   \cite{hac2}, Theorem 4.1).

\begin{proof}[Proof of Theorems \ref{main} and \ref{main2}]
We set up the notation in order to give a uniform presentation for all cases. Note first that by Proposition \ref{prop4}.a), under the hypotheses of Theorem \ref{main}, any divisor that represents $\ell$ is prime.

Let $L$ be an ample line bundle on $A$ that represents $\ell$, let $E$ be a prime divisor in $|L|$,  and let $D$ be a divisor in $m\ell$. We let $\cI_0=\cI_{Z_0}$ be the adjoint ideal $\cJ(A,E)$ and we let $\cI_1=\cI_{Z_1}$  be  the multiplier
ideal $\cI(A, \frac{1}{m} D)$. We have:
\begin{eqnarray*}
Z_0=\vide & \Longleftrightarrow &\hbox{the pair $(A,E)$ is canonical}\\
Z_1=\vide & \Longleftrightarrow &\hbox{the pair $(A, \frac1m D)$ is log terminal}
\end{eqnarray*}
so that we must prove (under suitable assumptions) that $Z_t$ is empty for $t\in\{0,1\}$.
We set as above,  for $P$ general in $\PA$,
$$ h_t=h^0(A,L\otimes  \cI_t\ot P)\in \{0,\dots,d\}$$
The point is to prove $h_t=d$ (Lemma \ref{lem1}.c)). 

We will use the following    reduced subvarieties of $\PA$ defined by  
$$V_i =\{ P\in \PA\mid H^i(A,L\otimes  \cI_t\otimes P)\ne 0\}
$$
which are analyzed in details in \S\ \ref{s2}, and set $V_{>0}  = \bigcup_{i>0}V_i  $.

\medskip
\noindent{\bf Case $t=0$.}  The exact sequence (\ref{adj}) shows that  Lemma \ref{bb} applies  with $\eps=1$ and $\cF=f_*\omega_X$, where
$f:X\to E$ is a desingularization.
 In particular, $V_{>0}\ne\PA$, hence 
 $$h_0= h^0(A,L\otimes  \cI_0\otimes P)=\chi(A,L\otimes  \cI_0\otimes P)=\chi(X,\omega_X)$$
 for $P$ general in $\PA$. Since $E$ is not fibered by (nonzero) abelian varieties, we obtain $h_0>0$ by\cite{EL}, Theorem 3.  

\medskip
\noindent{\bf Case $t=1$.} In this case,  $L\otimes\cI_1$ is a direct summand of the pushforward of a dualizing sheaf,\footnote{This can be seen as follows. Let $\mu:A'\to A$ be a log resolution of the pair $(A,D)$. Set $L'=\mu^*L\otimes \cO _{A'}(-\lfloor \frac {1}{m}\mu^*D\rfloor)$.
The divisor $$\mu^*D-m\Big\lfloor \frac {1}{m}\mu^*D\Big\rfloor\in |mL'|$$
defines a $\Z /m\Z$-cover $g:X\to A'$, and $X$ is normal with rational 
singularities. Let $\nu :X'\to X$ be a desingularization. One sees that $\omega _{A'}\otimes L'$ is a direct summand of 
$g_* \omega _X=g_*\nu_*\omega _{X'}$ (\cite{ev}, p. 33).

It follows that $\mu_*g_*\nu_*\omega _{X'}$ splits as a direct sum
of $m$ torsion free sheaves, one of these being  
\begin{eqnarray*}
\mu_* (\omega _{A'}\ot L')&=&\mu_*\Big(\omega _{A'/A}\ot \mu^*L\ot \cO _{A'}
\Big(-\Big\lfloor \frac {1}{m}\mu^*D\Big\rfloor\Big)\Big)\\
&=&L\ot \mu_* \Big(\omega _{A'/A}\ot\cO _{A'}
\Big(-\Big\lfloor \frac {1}{m}\mu^*D\Big\rfloor\Big)\Big)\\
&=&L\ot \cI_1 
\end{eqnarray*}} so that   Lemma \ref{bb} again applies,  with $\eps=0$. Moreover, since $g>d$ and  $\lfloor\frac1m D\rfloor=0$,   following the proof of \cite{hac1}, Theorem 1, one 
obtains $h_1>0$.

 \medskip
If $A$ is simple,  $V_{>0} $ is finite, and Lemma \ref{lem1}.e) implies Theorem \ref{main}.

We now prove Theorem \ref{main2}. Since $h_t>0$, we need only consider the case $d=2$.

If $h_t= 1$, since $V_1 $ has codimension at least $1$, the scheme $Z_t$ is a single point  and 
 $V_1  $ has dimension $\ge g-2$ (Lemma \ref{lem1}.d)). 
  This contradicts
the fact that $V_1 $ has dimension $0$
 (Lemma \ref{bb}.b)). 
 
 Hence $h_t=2$ and Theorem \ref{main2} is proved.
 \end{proof}

\medskip
We  now  interpret our results in terms of dimensions of loci of singularities.

\begin{theo}\label{th6}
 Let $(A,\ell)$ be a {\em simple} \pav\  of degree $d$
and dimension  $g> (d+1)^2/4$. Let $m$ and $k$ be positive integers. For all $D\in m\ell$,   we have
$$\dim \Sing_{mk}D< g-k$$
unless $k=1$ and $D=mE$, with $E\in \ell$.
\end{theo}

\begin{proof}According to Theorem \ref{main}, the hypotheses imply that the pair $(A,D)$ is canonical for $m=1$ and   the pair $(A, \frac{1}{m}D)$  is log terminal for $m\ge 2$, unless
$D=mE$, with $E\in \ell$. Since the pair $(A,E)$ is then also canonical, the theorem   follows from (\ref{sing}).
 \end{proof}

In the case of a polarization of degree $2$, we get a more precise result,  analogous to \cite{EL}, Corollary 2, and \cite{hac1},
Corollary 2.

\begin{theo}\label{th17}
 Let $(A,\ell)$ be an {\em indecomposable} \pav\  of degree $2$
and dimension $g>2$ and let $m$ and $k$ be positive integers. The
following properties are equivalent:
\begin{itemize}
\item[\rm (i)] for some $D$ in $ m\ell$, the locus
$\Sing_{mk}D$ contains an irreducible component of codimension  $k$ in
$A$;
\item[\rm (ii)] the \pav\ $(A,\ell)$ is a double \'etale cover of a product of $k$ nonzero \ppavs.
\end{itemize}
\end{theo}

\begin{proof}  If (ii) holds, the polarization is
represented by an
\'etale cover of the theta divisor of a product of $k$  nonzero   \ppavs,
hence (i) holds.

 Assume (i). By Theorem \ref{main2}.b) and (\ref{sing}), $D$ must have a
component of multiplicity at least $m$. If $D=mE$, with $E\in\ell$ prime, the pair $(A,E)$ is canonical (Theorem \ref{main2}.a)) and by (\ref{sing}), this contradicts (i). Therefore, 
by
Corollary \ref{ccc}.b),  there 
are  nonzero
\ppavs\ $(B_1,[\Theta_1])$ and $(B_2,[\Theta_2])$ and an isogeny  $p:A\to
B_1\times B_2$   such that $$D=mp^*(\Theta_1\times B_2)+ p^*( B_1\times D_2)
$$
 where $D_2\in |m \Theta_2|$.

 If $S$ is a component of $\Sing_{mk}D$ of maximal dimension, there is an integer $l\le k$ such that
$$S\subset \Sing_l \Theta_1\times \Sing_{m(k-l)}D_2
$$
From   \cite{kol}, Theorem 17.1,  we get 
 $$\codim (\Sing_l \Theta_1)\ge l$$
 From  \cite{EL}, Proposition 3.5,  we get 
 $$
\codim ( \Sing_{m(k-l)}D_2)\ge k-l$$
Since $S$ has codimension $k$, both inequalities must be equalities. By \cite{EL}, Corollary 2, $(B_1,[\Theta_1])$ splits as a product of $l$  nonzero   \ppavs, and by 
\cite{hac1}, Corollary 2,  $(B_2,[\Theta_2])$ splits as the product of at least $k-l$ \ppavs, so that   (ii) holds. 
\end{proof}

\section{Cohomological loci in $\PA$}\label{s2}

  Let $A$ be   an abelian variety. For any coherent sheaf $\cF$ on $A$ and integer $i$, we define reduced subvarieties of $\PA$ by setting
\begin{eqnarray*}
V_i(\cF)&=&\{ P\in \PA\mid H^i(A,\cF \otimes P)\ne 0\}\\
V_{>0}(\cF)&=&\bigcup_{i>0}V_i(\cF) 
\end{eqnarray*}
We investigate the geometry of these loci when $\cF$ is  the tensor product of an ample line bundle with an ideal sheaf, proving results that were used in the proof of Theorems \ref{main} and \ref{main2}.

 \begin{lemm}\label{lem1}Let $L$ be an ample line bundle of degree $d$ on  
an \av\ $A$ of dimension $g$, with base locus $\base|L|$, and let  $Z$ be a subscheme of $A$, with ideal sheaf $\cI$. Set $V_i= V_i(L\ot\cI)$ and
$$h=h^0 (A,L\otimes\cI\ot P) $$
for $P$ general in $\PA$. We have the following.
\begin{itemize}
\item[{\rm a)}] For $ i>\dim Z+1 $, the set $V_i$ is empty.
\item[{\rm b)}] If $h=0$, the set $V_{>0}$ is nonempty.
\item[{\rm c)}] We have $h\le d$, and $h=d$ \iff\  $Z$ is empty.
\item[{\rm d)}] If $h=d-1>0$, and if the \pav\ $(A,[L])$ is indecomposable,
 the
scheme $Z$ is finite and either $V_1=\PA$, or $Z$ is a single (reduced) point
$z$ and $V_1=\phi_L(\base|L|-z)$, so that $\dim V_1\ge g-d$.
\item[{\rm e)}]   If $A$ is simple and $0<h<d$,  we have $\dim Z\le d-1-h$ and  $\dim V_{>0}\ge g-(d+1)^2/4$.
\end{itemize}
\end{lemm}

  Note that
$h=h^0 (A,L\otimes \cI\ot P)=\chi (A,L\otimes\cI)  $ for $P\notin V_{>0} $.

\begin{proof} For $ i>\dim Z+1 $, we have $H^i(A,L\ot\cI \otimes P)\isom H^i(A,L \otimes P)=0$ for all $P\in\PA$. This proves a).

If $V_{>0}$ is empty,   
$$h=h^0 (A,L\otimes \cI\ot P)=\chi (A,L\otimes\cI)  $$
 for all $P\in\PA$. If $h=0$, we have $H^i (A,L\otimes \cI\ot P)=0$
 for all integers $i$ and all $P\in\PA$. This is impossible by \cite{muk}, Corollary 2.4,   and b) is proved.

We have $h= d$  \iff\   all
sections of $L$ vanish on general translates of $Z$; this happens \iff\
$Z$ is empty. This   proves c).

Let us now prove d). Set
$$J=\{ (s,a)\in \P H^0(A,L)\times A\mid s\vert_{Z+a}\equiv 0\}$$
 The fiber of a point $a$ of $A$ for the   second
projection $q:J\to A$ is isomorphic to
$\P H^0(A,L\otimes P_{\phi_L(a)}\otimes\cI)$. If $h>0$,  a unique irreducible component $I$ of $J$ dominates $A$, and  
 $\dim I=  g+h-1$. 

 Let $p: I\to \P H^0(A,L)$ be the first projection. Since any nonempty  
$F_s=q\bigl( p^{-1}(s)\bigr)$   satisfies
$Z+F_s\subset \div (s)$, we have 
$$g-1\ge \dim F_s \ge \dim I -\dim p(I) \ge g+h-1-\dim p(I)\ge g-(d-h)
$$

Assume $h=d-1>0$. Then $p$ is surjective and $F_s$ has
dimension $g-1$. The divisor of a general section
$s$ being prime, the inclusion $z+F_s\subset
\div (s)$ is an equality for all $z$ in $Z$. This
implies that $Z$ is finite. If $V_1\ne\PA$, the length of $Z$ is $d-h=1$, so that $Z=\{z\}$, and an element $P=P_{\phi_L(a)}$ of $\PA$ satisfies  
$H^1(A,L\otimes P \otimes\cI)\ne0$ \iff\ the restriction 
$$H^0(A,L\otimes P )\to H^0(Z,L\otimes \cO_Z\ot P)\isom \C_z$$ is not
surjective; in other words, if all sections of $L\otimes P$ vanish at $z$, \ie,
$z\in\base|L\otimes P|=\base|L|-a$. This   proves d).

Assume $A$ is {\em simple} and $0<h<d$. Then $Z$ is nonempty and the inclusion $Z+F_s\subset
\div (s)$  implies (\cite{De3}, Corollaire 2.7)
$$\dim Z\le g-1-\dim F_s\le d-1-h$$
For $a$ general in $A$, the subvariety $p(q^{-1}(a))$ of $ \P H^0(A,L) $ is a linear subspace of dimension $h-1$. It must vary with $a$, because a nonzero $s$ does not vanish on all translates of $Z$. It follows that the linear span of $p(I)$ has dimension at least $h$.
For $s_1,\dots,s_{h+1}$ general elements in $p(I)$,
one has (\cite{De3}, Corollaire 2.4)
$$\dim (F_{s_1}\cap \dots \cap F_{s_{h+1}})\ge  g-(h+1)(d-h)\ge g-(d+1)^2/4$$
For $a\in F_{s_1}\cap \dots \cap F_{s_{h+1}}$, the sections  $s_1,\dots,s_{h+1}$ all vanish on $Z+a$, hence  $h^0(A,L\otimes \cI \ot P _{\phi _L(a)})\ge h+1$. This implies $P _{\phi _L(a)} \in V_{>0} $  and   proves e).
\end{proof}

Assume now that there is a smooth variety $X$ with a morphism $f:X\to A$
such that the sheaf $\cF$ on $A$ is a direct summand of $f_*\omega_X$.
Let $B$ be an abelian variety with a morphism $\pi:A\to B$. For all integers $i$ and $j$, and any torsion point $P_0\in\PA$, {\em  every irreducible component of  $V_i(R^j\pi_*(\cF\ot P_0))$  is an abelian subvariety of $\PB$ of codimension at least $i$ translated by a torsion point.}
This applies in particular to the loci $V_i(\cF)$ in $\PA$. 

When $P_0=0$, this is a particular case of \cite{hp}, Theorem 2.2. For the general case, associate to the torsion element $f^*P_0\in \Pic ^0(X)$
a cyclic \'etale cover   $p:X'\to X$. Then $\omega_X\ot f^*P_0$ is a direct summand of $p_*\omega_{X'}$, hence $\cF\ot P_0$  is a direct summand of $f_*p_*\omega_{X'}$.

\begin{lemm}\label{bb}
Under the hypotheses and notation of Lemma \ref{lem1}, assume further that there is
an exact sequence
\begin{equation*}
0\to \cO_A^{\oplus \eps}\to L\ot\cI\to \cF\to 0
\end{equation*}
for some $\eps\in\N$, where $\cF$ is a direct summand of a pushforward of a dualizing sheaf.  The following properties hold.
\begin{itemize}
\item[{\rm a)}]Every irreducible component of  $V_i(L\ot \cI )$  is an
abelian subvariety of $\PA$ of codimension at least $i$ translated  by a torsion point.
\item[{\rm b)}] If the support of $\cF$ is not contained in any nonample divisor of $A$, we have, for any  $i>0$ such that $V_i(L\ot \cI ) $ is nonempty,  $\dim Z\ge i-1+\dim V_i(L\ot \cI ) $.
\end{itemize}
\end{lemm}

\begin{proof}Since
$V_i(L\otimes \cI)\moins\{0\}=V_i(\cF)\moins\{0\}$, item
 a) holds. Let us prove b).
Since b) follows from Lemma \ref{lem1}.a) when $V_i(L\ot \cI ) $ is finite, we may pick   a common irreducible component $V$ of 
$ V_i(\cF)$ and of $V_i(L\ot \cI ) $ of
maximal {\em positive} dimension. Let $B= \PV$ and let $\pi :A\to B$
be the induced morphism. Let $P_0\in\PA$ be a torsion point such that  $V=P_0+\pi^*\PB$.

We know that $V_k(R^j\pi_*(\cF\otimes P_0))$ has codimension at least $k$ in $\PB$. It follows that for
general $P \in \PB$, and all $k>0$ and $ j\geq 0$,
$$H^k(B,R^j\pi _* (\cF \ot P_0)\ot P )=0$$
hence
\begin{equation}\label{rr}
H^0( B,R^i\pi _* (\cF \ot P_0)\ot P )\isom H^i(A,\cF\ot P_0 \ot
\pi ^*P )\ne 0
\end{equation}
because $P_0 \ot\pi ^*P$ is in $V\subset V_i(\cF )$. We have an exact sequence
$$R^i\pi_* (P_0 )^{\oplus \eps}\lra R^i\pi _* (L\ot\cI \ot P_0) \lra R^i\pi _* (\cF \ot P_0)\stackrel{\delta}{\lra} R^{i+1}\pi_* (P_0 )^{\oplus \eps}$$
The sheaves $ R^j\pi_* (P_0 )$ are  direct sums of numerically trivial line bundles on $B$ (this follows from the proof of \cite{Kempf},  Theorem 1). By a result of Koll\'ar (\cite{kol3}, Theorem 3.4; \cite{hp}, Theorem 2.1), the sheaf 
$R^i\pi _* (\cF \ot P_0)$  is torsion-free on  $\pi(\Supp \cF)$, which is $B$ by hypothesis. It follows that the support of the sheaf $R^i\pi _* (L\ot\cI \ot P_0)$ is $B$: if it is not, the map $\delta$ is generically injective, hence injective; but a twist of $ R^{i+1}\pi_* (P_0 )^{\oplus \eps}$ by a general $P\in\PB$ has no nonzero section, contradicting (\ref{rr}).

Since $R^i\pi_*(L\ot P_0)=0$, the short exact sequence
\begin{equation*}
0\to L\ot \cI  \ot P_0\to L\ot P_0\to L\otimes\cO_Z\ot P_0\to 0
\end{equation*}
yields a surjection
\begin{equation*}
 R^{i-1}\pi_*(L\otimes\cO_Z\ot P_0)\twoheadrightarrow
R^i\pi_*(L\ot \cI  \ot P_0) 
\end{equation*}
 hence the support of $R^{i-1}\pi_*(L\ot \cO_Z \ot P_0)$ is also $B$. In particular, all fibers of $\pi\vert_Z:Z\to B$ have dimension at least $i-1$, and b) follows.
\end{proof}

\section{Reducible divisors in indecomposable polarizations}\label{app}

We gather in this last section elementary results on \pavs\  that were used earlier.

\begin{lemm}\label{lem77}
 Let $(A,\ell)$ be an indecomposable  \pav. If the restriction of $\ell$
  to an abelian subvariety  
$A_1$   of $A$ is principal, either $A_1=0$ or $A_1=A$.
\end{lemm}

\begin{proof}Let $A_2$ be the neutral component of the kernel of
$$f: A\stackrel{\phi_\ell}{\lra} \PA\lra\Pic^0(A_1)
$$
The sum map $p:A_1\times A_2\to A$ is an isogeny of \pavs\ whose kernel is 
isomorphic to $A_1\cap A_2$. Since
$\ell$ induces a principal polarization on $A_1$, the restriction of
$f$ to $A_1$ is injective, \ie, $A_1\cap A_2=\{0\}$ and $p$ is injective.\end{proof}

We now study indecomposable polarizations that contain nonprime divisors.  The situation is manageable when $A$ is simple or the degree is $2$. For degrees at least $3$, more and more exceptional cases arise.

\begin{prop}\label{prop4}
 Let $(A,\ell)$ be an indecomposable \pav\  of degree
$d$ and dimension
$g\ge d$ such that $\ell$ contains a nonprime divisor $E$. 

{\rm a)} The \av\ $A$ is {\em not} simple.

{\rm b)} If $d=2$,  there exist a {\em decomposable}
\ppav\ $(B,[\Theta])$   and an isogeny  $p:A\to
B$    of degree
$d$  such that $E=p^*\Theta$.
\end{prop}

\begin{proof}
Write
$E=E_1+E_2 $, with $E_1$ and $E_2$   effective and nonzero, and let, for $j\in\{1,2\}$,
$A_j= K(E_j)^0$,
$B_j=A/A_j$, and  $g_j=\dim B_j>0$.   Since $A_1\cap A_2$ is contained in $K(\ell)$, it is finite,  hence $g_1+g_2\ge g$. 
 There is an ample divisor
$D_j$ on
$B_j$ which pulls back to
$E_j$, and
$$d=\frac{1}{g!}(E_1+E_2)^g=\frac{1}{g_1!(g-g_1)!}E_1^{g_1}E_2^{g-g_1}+\cdots+
\frac{1}{(g-g_2)!g_2!}E_1^{g-g_2}E_2^{g_2}
$$
 The first term of this sum is 
\begin{eqnarray*}
\frac{1}{(g-g_1)!}\deg(D_1)[A_1]\cdot E_2^{g-g_1}&=&\deg(D_1)\deg(E_2\vert_{A_1})\\
&=&
\deg(D_1)\deg(\ell\vert_{A_1})>0
\end{eqnarray*}
and similarly for the last term, hence
all terms  are
positive\footnote{This follows for example from the Teissier--Hovanski
inequalities
$$E_1^i E_2^{g-i}\ge
\bigl(E_1^{g_1}E_2^{g-g_1}\bigr)^{\frac{i-g+g_2}{g_1+g_2-g}}
\cdot\bigl(E_1^{g-g_2}E_2^{g_2}\bigr)^{\frac{g_1-i}{g_1+g_2-g}}
$$
for $g_1\ge i\ge g-g_2$.} integers.

When $A$ is simple, we have $g_1=g_2=g$, hence $d>g$. This proves a).

We now assume $d=2$ and  prove b).
By Lemma \ref{lem77}, $\ell$ does not  restrict to a principal
polarization on  $A_j$ unless $A_j=0$, \ie, $g_j=g$.
The only possibility is  $g_1+g_2=g$, the polarization  $[D_j]$
 on $B_j$ is principal, the map $p:A\to B_1\times B_2$ is
an isogeny,   and $E=p^*(D_1\times B_2)+ p^*(B_1\times D_2) $.
 \end{proof}

We use these results to bound the multiplicities of the components of elements of $m\ell$.

\begin{coro}\label{ccc}
 Let $(A,\ell)$ be an indecomposable \pav\  of degree
$d$ and dimension
$g\ge d$ and let $D\in m\ell$, with $m>0$.
\begin{itemize}
\item[{\rm a)}] If $A$ is   simple, we have $\lfloor\frac{1}{m+1}D\rfloor=0$. Moreover, $\lfloor\frac{1}{m}D\rfloor=0$ unless $D=mE$, with $E\in\ell$.
\item[{\rm b)}] If $d=2$, we have $\lfloor\frac{1}{m+1}D\rfloor=0$. Moreover, $\lfloor\frac{1}{m}D\rfloor=0$  unless $D=mE$, with $E\in\ell$ prime, or  there 
are  nonzero
\ppavs\ $(B_1,[\Theta_1])$ and $(B_2,[\Theta_2])$,  and an
isogeny 
$p:A\to B_1\times B_2$    of degree
$2$,  such that
\begin{equation*}
D=mp^*(\Theta_1\times B_2)+ p^*( B_1\times D_2)
\end{equation*}
 with $D_2\in |m \Theta_2|$. 
 \end{itemize}
\end{coro}

\begin{proof}  The arguments of the proofs of Lemmas 2.2 and 2.3 of \cite{hac1} yield:
\begin{itemize}
\item  $\lfloor\frac{1}{m+1}D\rfloor=0$,
unless  $\ell$ contains a reducible divisor of the form $E_1+E_2$, with $E_2$ ample and $E_1$,  $E_2$, and $D-(m+1)E_1$ effective nonzero;
\item  $\lfloor\frac{1}{m}D\rfloor=0$,
unless   $\ell$ contains a divisor of the form $E_1+E_2$, with $E_1$,  $E_2$, and $D-mE_1$ effective.
\end{itemize}
The corollary therefore follows from the proposition.
\end{proof}

 \begin{rema}\upshape\label{rem10}
Polarized abelian varieties that satisfy the condition in Proposition \ref{prop4}.b) are all obtained as follows. Let $d$ be any positive integer and, for $j\in\{1,2\}$, let $(B_j,[\Theta_j])$ be a nonzero \ppav\ of
dimension 
$g_j$. Endow $B= B_1\times
B_2$ with the product polarization $[\Theta]$.
Choose  a point $\b_j$ in $ B_j$ 
of order
$d$  and consider the cyclic isogeny $p:A\to B$ of degree $d$
associated with the point
$\phi_{\Theta}(\b_1,\b_2)$ of
$\Pic^0(B)$. The divisor $ p^*\Theta$ is reducible and defines
a polarization $\ell$  of type
$(d)$ on $A$.

Let $p_j:A_j\to B_j$  be the degree $d$ cyclic isogeny associated with $\phi_{\Theta_j}(\b_j)$   and let $\ell_j$ be the polarization $[p_j^*\Theta_j]$ (of type $(d)$) induced  on $A_j$. There is a
factorization
$$p_1\times p_2:A_1\times A_2\stackrel{\pi}{\lra} A
\stackrel{p}{\lra} B
$$
Another way to construct $(A,\ell)$ is to start from nonzero \pavs\ $(A_1,\ell_1)$ and $(A_2,\ell_2)$ of type $(d)$, to choose elements
$\a_1\in K(\ell_1)$ and $\a_2\in K(\ell_2)$ or order $d$, and to take the quotient of $A_1\times A_2$ by the subgroup generated by  
$(\a_1,\a_2)$. For more details, see \cite{Deduke}, Proposition 9.1.

Assume $d$ is {\em prime.} A \pav\ of degree $d$ is decomposable \iff\ it has a nonzero principally polarized abelian factor. It follows that the \pav\ 
 $(A,\ell)$ obtained by the above construction is   indecomposable \iff\
both
\pavs\ $(A_1,\ell_1)$ and $(A_2,\ell_2)$ are indecomposable.
 \end{rema}

\end{document}